\documentclass{article}
\usepackage{amsmath,amssymb,amsthm,amscd, mathtools, latexsym}
\usepackage{enumerate,varioref, dsfont, fancyhdr}
\usepackage{tikz-cd}
\usepackage{enumitem}
\usepackage{multicol}
\usepackage{hyperref}
\usepackage{url}
\usepackage{textcomp}

\newtheorem{Thm}{Theorem}[section]

\newtheorem{Lem}[Thm]{Lemma}
\newtheorem{Cor}[Thm]{Corollary}
\newtheorem{Conj}[Thm]{Conjecture}
\newtheorem{Ques}[Thm]{Question}
\newtheorem{Ass}[Thm]{Assumption}

\newtheorem{Prob}[Thm]{Problem}
\newtheorem{Def}[Thm]{Definition}
\newtheorem{Rem}[Thm]{Remark}
\newtheorem{Not}[Thm]{Notation}

\theoremstyle{definition} 

\def\cit{{\mathbb C}}

\def\qit{{\mathbb Q}}
\def\zit{{\mathbb Z}}

\def\pit{{\mathbb P}}

\def\0{{\mathcal O}}

\def\Hom{\mathop{\rm Hom}\nolimits}

\def\X{{\mathcal X}}

\usepackage[title]{appendix}

\begin{document}
\title{On the unramified cohomology of certain quotient varieties}
\author{H. Diaz}
\newcommand{\Addresses}{{\bigskip \footnotesize
\textsc{Department of Mathematics, University of California Riverside, Riverside, CA 92521} \par \nopagebreak
\textit{Email address}: \ \texttt{humbertoadiaziii@gmail.com}}}

\date{}
\maketitle

\begin{abstract}
\noindent In this note, we consider unramified cohomology with $\zit/2$ coefficients for some (degree two) quotient varieties and describe a method that allows one to prove the non-vanishing of these groups under certain conditions. We apply this method to prove a non-vanishing statement in the case of Kummer varieties. Combining this with work of Colliot-Th\'el\`ene and Voisin, we obtain a new type of three-dimensional counterexample to the integral Hodge conjecture.
 
\end{abstract}
\date{}

\noindent Let $X$ be a smooth projective variety over $\cit$. We will consider the {\em unramified cohomology with coefficients in $A$}, $H^{i}_{nr} (X, A) = \Gamma(X, \mathcal{H}_{X}^{i}(A))$, where $\mathcal{H}_{X}^{i}(A)$ denotes the Zariski sheaf over $X$ associated to the presheaf $U \mapsto H^{i} (X, A)$. When $A = \zit/n$, it has been known since \cite{CTO} that these groups are (stably) birational invariants, and this been used quite often (and by many different authors) to prove the existence of unirational varieties are not (stably) rational.\\
\indent Another application of these groups is to the integral Hodge conjecture, which asserts that every $\alpha \in H^{p,p} (X, \zit)$ is algebraic. When $p= 0, dim(X)$, this is trivially true and when $p=1$ this is true by the Lefschetz $(1,1)$-theorem. For all other $p$, it has been known since the counterexamples of Atiyah and Hirzebruch \cite{AH} that the conjecture is false. Of particular interest is the case when $X$ is a threefold (necessarily for $p=2$). In this direction, various results have been produced, both positive and negative. On the positive side, Voisin \cite{V} proved that the integral Hodge conjecture holds for uniruled threefolds and for Calabi-Yau threefolds. Grabowski \cite{Gr} also proved the conjecture holds for Abelian threefolds. On the negative side, Koll\'ar \cite{Ko} produced the first counterexamples for threefolds, in the form of non-algebraic $(2,2)$ cohomology classes on general hypersurfaces in $\pit^{4}$ of sufficiently large degree. Other counterexamples have been found with Kodaira dimension $1$ (Colliot-Th\'el\`ene and Voisin \cite{CTV}; Totaro \cite{T}). More recently, Benoist and Ottem \cite{BO} produced threefolds of Kodaira dimension zero that fail the Hodge conjecture, the counterexamples they give being products of Enriques surfaces with very general elliptic curves. Their counterexample was generalized by Shen \cite{Sh} to products of Enriques surfaces with very general odd degree hypersurfaces of higher dimension. Colliot-Th\'el\`ene \cite{CT} reinterpreted the result of Benoist and Ottem using unramified cohomology (together with a degeneration technique of Gabber \cite{CTG}) and gave further counterexamples in dimension $3$. This relation between the integral Hodge conjecture and unramified cohomology was first explored in \cite{CTV}, and the main results (in principle) give a recipe for producing counterexamples to the integral Hodge conjecture if one can exhibit non-trivial unramified cohomology classes (in degree $3$).\\
\indent Unfortunately, the drawback to working with unramified cohomology groups is that they are quite mysterious in degree $\geq 3$. Indeed, the canonical map
\begin{equation} H^{i} (X, \zit/n) \to H^{i}_{nr} (X, \zit/n)\label{canonical} \end{equation}
is surjective for $i \leq 2$. By contrast, it is a non-trivial question in general whether or not the map (\ref{canonical}) is even non-vanishing. One technique that works to prove the non-vanishing of (\ref{canonical}) (under suitable conditions) was developed by Bloch and Esnault \cite{BE}; this is a mixed characteristic approach that is typically used in proving results about non-divisibility of algebraic cycles. As we discuss below, it implies that when $X$ is a very general Abelian variety of dimension $g$, the map (\ref{canonical}) is non-vanishing for $1 \leq i \leq g$.\\
\indent The use of the Bloch and Esnault method until now has been restricted to instances for which $H^{0,i} (X) \neq 0$. In this note, we would like to develop a strategy involving this method that works to prove the non-vanishing of $H^{i}_{nr} (X, \zit/n)$ in some instances for which $H^{0,i} (X) = 0$. We will be concerned with the case $n=2$ and that $X$ is a degree two quotient variety. While our main application is to prove non-vanishing statements about $H^{i}_{nr} (X, \zit/2)$ when $X$ is a Kummer variety (see Corollary \ref{kum-cor}), the strategy developed here is flexible and can be applied to other types of degree two quotients. One consequence of these computations (in conjunction with the results of \cite{CTV}) is what is (evidently) a new type of three-dimensional counterexample to the integral Hodge conjecture:
\begin{Thm}\label{main} There exists a smooth projective simply-connected threefold $X$ of Kodaira dimension zero violating the integral Hodge conjecture. In fact, one can find such an $X$ defined over $\qit$.
\end{Thm} 
\noindent The counterexample given is a Kummer variety $X$ of dimension $3$ and is curious for a few reasons. First, by an old result of Spanier \cite{Sp} (discussed below), not only is $X$ simply connected, but $H^{odd} (X, \zit) =0$ and $H^{*} (X, \zit)$ is torsion-free (so that, in particular, this is a non-torsion counterexample). So, we will obtain that (for many such $X$) $H^{3}_{nr} (X, \zit/2(2)) \neq 0$ even though $H^{3} (X, \zit/2(2)) = 0$. To get around the fact that $H^{3} (X, \zit/2(2)) = 0$, we show that $H^{3}_{nr} (X, \zit/2(2)) = 0$ by a rather na\"ive ``descend-and-extend" argument (see \S 2): using the corresponding Abelian variety $A$, we descend unramified cycles to a suitable open subset of $X$ and then show that cycles obtained in this way must extend to unramified cycles on $X$. \\
\indent Additionally, as we note above, the integral Hodge conjecture holds for many threefolds with trivial canonical class. The counterexample obtained in \cite{BO} is of the form $S \times E$, where $S$ is an Enriques surface and $E$ is a very general elliptic curve, so that the canonical divisor is $2$-torsion. This phenomenon is again evident in the case of a Kummer threefold, albeit in different way; i.e., the minimal model (the singular Kummer threefold) has $2$-torsion canonical divisor. \\
\indent Finally, as noted above, the counterexamples we obtain can be defined over number fields. This is one advantage of using the Bloch and Esnault method in that it applies to varieties that are very general in moduli, as well to those defined over number fields. The problem of finding $3$-dimensional counterexamples to the integral Hodge conjecture defined over number fields was carefully considered by Totaro in \cite{T} who used the Hassett-Tschinkel method (among others) to obtain many such counterexamples. Degeneration arguments can also be used to obtain counterexamples defined over number fields (as in op. cit. and \cite{BO} Rem. 2.2). \\
\indent It is natural to speculate if the above technique will also work to answer (positively) the following much more elusive question:
\begin{Ques} Does there exist a simply-connected threefold with a non-algebraic torsion cohomology class?
\end{Ques}
\noindent Our plan will be as follows. In section $1$, we review some essential properties of unramified cohomology and the method of Bloch and Esnault. In section $2$, we describe the ``descend-and-extend" strategy mentioned above and its applications. In section $3$, we prove Theorem \ref{main}. In section $4$, we address how Theorem \ref{main} (conjecturally) generalizes to higher dimension.
\subsection*{Acknowledgements} The author would like to thank Olivier Benoist for his feedback and for finding a mistake in a previous draft (that led to this draft). He would also like to thank John Christian Ottem for his interest.
\subsection*{Notation} Unless otherwise specified, all varieties (reduced schemes of finite type) considered below are defined over $\cit$. Moreover, we let the cohomology group $H^{*} (-, A)$ denote singular cohomology with $A$ coefficients. As usual, we let $\zit(1) = \zit(2\pi i) \subset \cit$ and, for $m \geq 0$, let $\zit(m) = \zit(1)^{\otimes m}$. Additionally, for any integer $m$ and Abelian group $A$, we define
\[ A(m) = \left\{\begin{array}{ll} A \otimes \zit(m) & \text{ for } m \geq 0\\
\Hom(\zit(m), A) & \text{ for } m < 0 \end{array} \right.\]  
With this notation, we have a canonical isomorphism $\zit/n(m) \cong \mu_{n}^{\otimes m}$, where $\mu_{n}$ denotes the group of $n^{th}$ roots of unity and $m\geq 0$. Finally, for $A$ an Abelian group, we denote its $n$-torsion subgroup by $A[n]$.
\section{Preliminaries}

\begin{Def} Let $X$ be a smooth variety over $\cit$, $A$ be an Abelian group and $\mathcal{H}_{X}^{i}(A)$ be the Zariski sheaf over $X$ associated to the presheaf $U \mapsto H^{i} (U, A)$. Then,
\[ H^{i}_{nr} (X, A) := H^{0} (X, \mathcal{H}_{X}^{i}(A)) \]
is the {\em $i^{th}$ unramified cohomology group with coefficients in $A$}.
\end{Def}
\noindent When $A$ is torsion, we have the following alternative characterization. Indeed, note that Theorem 4.2 of \cite{BlO} gives the following well-known short exact sequence:
\begin{equation}  H^{i}_{nr} (X, A) \hookrightarrow \bigoplus_{x \in X^{(0)}} H^{i} (\cit(x), A) \to \bigoplus_{x \in X^{(1)}} H^{i-1} (\cit(x), A(-1)) \label{res}\end{equation}
Here, $X^{(n)}$ denotes the set of irreducible codimension $n$ subvarieties on $X$ and 
\begin{equation} H^{*} (\cit(x), A) := \mathop{\lim_{\longrightarrow}}_{U \subset x} H^{*} (U, A)\label{Galois}\end{equation}
in which the limit ranges over all Zariski open subsets of $x$ (viewed as a subvariety of $X$) and where the rightmost arrow of (\ref{res}) is the sum of residue maps. This somewhat ad hoc definition (\ref{Galois}) is used in loc. cit. as a suitable generalization of Galois cohomology with torsion coefficients. An immediate consequence of this characterization is that the restriction map $H^{i}_{nr} (X, A) \to  H^{i}_{nr} (U, A)$ is an isomorphism when $U \subset X$ is an open subset whose complement has codimension $\geq 2$ in $X$. Moreover, if the complement of $U$ is a closed subvariety $D$ of codimension $1$, there is a short exact sequence:
\begin{equation} 0 \to H^{i}_{nr} (X, A) \to  H^{i}_{nr} (U, A) \to \bigoplus_{x \in D^{(0)}} H^{i-1} (\cit(x), A(-1)) \label{kernel}\end{equation}
where the right non-zero arrow is the sum of residue maps. Moreover, there is the natural map:
\begin{equation}  H^{i} (X, A) \to H^{i}_{nr} (X, A)\label{obvious}\end{equation}
When $A = \zit/n(m)$, there is the local-to-global spectral sequence:
\begin{equation} E_{2}^{p,q} = H^{p}_{Zar} (X, \mathcal{H}^{q}_{X} (\zit/n(m))) \Rightarrow H^{p+q} (X, \zit/n(m))\label{spectral} \end{equation}
for which $E_{2}^{p,q} = 0$ for $p>q$ thanks to the Gersten resolution obtained in \cite{BlO}. Moreover, the natural map \ref{obvious} arises as an edge map that is an isomorphism for $i \leq 1$, and in general there is a short exact sequence:
\begin{equation} 0 \to N^{1}H^{i} (X, \zit/n(m)) \to H^{i} (X, \zit/n(m)) \to H^{i}_{nr} (X, \zit/n(m))\label{niveau}\end{equation}
here $N^{*}$ denotes the coniveau filtration on $H^{*} (-, \zit/n(m))$. When $i=2$ the right arrow is surjective so that, in particular, it follows that 
\begin{equation} H^{2}_{nr} (X, \zit/n(1)) \cong Br(X)[n] \label{Brauer}\end{equation}
using the Kummer exact sequence. For $i \geq 3$, the coniveau filtration is quite mysterious and it is unclear in general whether or not 
\[N^{1}H^{i} (X, \zit/n(m)) \neq H^{i} (X, \zit/n(m))\] 
(and, hence, whether or not (\ref{obvious}) vanishes) even when in instances for which $N^{1}H^{i} (X, \zit(m)) \neq H^{i} (X, \zit(m))$. To prove non-vanishing of (\ref{obvious}), a few techniques exist. There is the degeneration method given in the appendix of \cite{CTG} and used in \cite{CT}. The idea is that if one is able to spread out $\gamma \in H^{i} (X, \zit/n(m))$ to a cycle $\Gamma \in  H^{i} (\mathcal{X}, \zit/n(m))$, where $\X \to S$ is a family of smooth irreducible varieties containing $X$ as a (very general) fiber, then the set of $s \in S$ for which 
\begin{equation} \Gamma\in  N^{1} H^{i} (\mathcal{X}_{s}, \zit/n(m)) \label{Gamma} \end{equation} is a $G_{\delta}$ subset of $S$. In particular, if one can find some $s \in S$ for which (\ref{Gamma}) does not hold, then a Baire category argument shows that very generally (\ref{Gamma}) does not hold. This implies the non-vanishing of (\ref{obvious}), since $X$ was very general.\\
\indent There is also the following method of Bloch and Esnault. As noted in the introduction, this is a mixed characteristic method and has the advantage that it does not require spreading out and, hence, can be applied to prove the non-vanishing of (\ref{obvious}) for varieties over number fields. Its proof uses the spectral sequence in \cite{BK} for $p$-adic \'etale cohomology that degenerates for smooth projective varieties over a $p$-adic field with good ordinary reduction. This result has been used repeatedly in the context of finding non-divisible cycles in the Griffiths group. 
\begin{Thm}[Bloch-Esnault, \cite{BE} Theorem 1.2]\label{BE} Let $X$ be a smooth projective irreducible variety over a complete discrete valuation field $K$ with perfect residue field $k$ of mixed characteristic $(0,p)$. Suppose further that $X$ has good ordinary reduction and that 
\begin{enumerate}[label=(\alph*)]
\item\label{no-tors-2} The crystalline cohomology of the special fiber $Y$ has no torsion.
\item\label{forms} $H^{0} (Y, \Omega^{m}_{Y}) \neq 0$.
\end{enumerate}
Then, $N^{1}H^{m}_{\text{\'et}} (X_{\overline{K}}, \zit/p) \neq H^{m}_{\text{\'et}} (X_{\overline{K}}, \zit/p)$.
\end{Thm}
\noindent As an immediate application, we have the following consequence for Abelian varieties over number fields.
\begin{Cor}\label{Cor-to-BE} Suppose that $A$ is a complex Abelian variety that is defined over a number field which has good ordinary reduction at some prime dividing the prime $p \in \zit$. Then, for all $1 \leq i \leq g$, the map $H^{i} (A, \zit/p) \to H^{i}_{nr} (A, \zit/p)$ is non-zero.
\end{Cor}
\begin{Rem}\label{many} As noted on p. 108 of \cite{BK}, Deligne proved that the set of ordinary hypersurfaces in projective space (over a $p$-adic field) of any given degree make up an open dense set in the moduli space. This implies that over a given number field $k$, there is a Zariski dense subset of the moduli of plane curves defined over $k$ that have good ordinary reduction at a prime dividing $p$. In particular, this is true for the set of plane curves over $k$ of degree $4$, for which the corresponding set of Jacobians forms a Zariski open subset of the moduli of Abelian threefolds. We deduce that over a given number field $k$, there is a Zariski dense subset of the moduli of Abelian threefolds defined over $k$ that have good ordinary reduction at a prime dividing $p$. (We will be interested in the case that $p=2$.)
\end{Rem}
\section{Main strategy}

\subsection{Descent and extension}
\begin{Not} Suppose that $U$ is a smooth variety over $\cit$ and let $\pi: V \to U$ be a degree $2$ finite \'etale cover. Denote the corresponding Galois group by $C_{2} = Aut(Y/X)$ and let $\iota$ be the involution that generates it. 
\end{Not}
\noindent Our first lemma is the descent part of the strategy. It is easy and shows that mod $2$ unramified cohomology cycles descend along a double cover in degrees for which the pull-back map is surjective. 
\begin{Lem}\label{litmus} Suppose that the natural map 
\[ H^{i} (V, \zit/2) \to H^{i}_{nr} (V, \zit/2) \] 
does not vanish and that the pull-back $H^{i} (U, \zit/2) \xrightarrow{\pi^{*}} H^{i} (V, \zit/2)$ is surjective. Then, the natural map $H^{i} (U, \zit/2) \to H^{i}_{nr} (U, \zit/2)$ does not vanish.
\begin{proof}  Consider the obvious commutative diagram:
\[\begin{tikzcd}
H^{i} (U, \zit/2) \arrow{r}{\pi^{*}} \arrow{d} & H^{i} (V, \zit/2)\arrow{d} \\
H^{i}_{nr} (U, \zit/2) \arrow{r}{\pi^{*}}  & H^{i}_{nr} (V, \zit/2)
\end{tikzcd}
\]
Since the right vertical arrow is non-vanishing and the top horizontal arrow is surjective, it follows that the left vertical arrow is also non-vanishing, as desired.
\end{proof}
\end{Lem}
\noindent If one assumes the surjectivity of $H^{i} (U, \zit/2) \xrightarrow{\pi^{*}} H^{i} (V, \zit/2)$, a necessary consequence is that $\iota$ act trivially on $H^{i} (V, \zit/2)$. This triviality of $\iota^{*}$ is not really much of a restriction; for instance, any involution that acts diagonalizably on $H^{i} (V, \zit)$ (with $H^{i+1} (V, \zit)$ is torsion-free) acts trivially on $H^{i} (V, \zit/2)$. One can ask if this is sufficient to ensure the surjectivity assumption in Lemma \ref{litmus}:
\begin{Prob}[\'Etale descent]\label{hard} Suppose that $\iota$ acts trivially on $H^{i} (V, \zit/2)$. Is the pull-back map $H^{i} (U, \zit/2) \xrightarrow{\pi^{*}} H^{i} (V, \zit/2)$ surjective?
\end{Prob}
\noindent Besides the trivial case that $i=0$, it turns out the answer is conditionally ``yes" for $i=1$, provided that we impose the following extra assumption:
\begin{Ass}\label{ass} The natural short exact sequence:
\begin{equation} 1 \to \pi_{1} (V) \to \pi_{1} (U) \to C_{2} \to 1\label{want-split} \end{equation}
splits on the right; i.e., $\pi_{1} (U) \cong \pi_{1} (V) \rtimes C_{2}$.
\end{Ass}
\begin{Lem}\label{lem-1} If $\iota$ acts trivially on $H^{1} (V, \zit/2)$ and Assumption \ref{ass} holds, then the pull-back map \[ H^{1} (U, \zit/2) \xrightarrow{\pi^{*}} H^{1} (V, \zit/2)\] 
is surjective.
\begin{proof} Since $\pi_{1} (U) = \pi_{1}(V) \rtimes  C_{2}$ by assumption, we need to show that the restriction map in group cohomology:
\[ H^{1} (\pi_{1}(V) \rtimes C_{2}, \zit/2) \to H^{1} (\pi_{1}(V), \zit/2) \]
is surjective. To this end, since $\iota$ acts trivially on $H^{1} (V, \zit/2)$ by assumption, inflation-restriction gives an exact sequence:
\[ H^{1} (\pi_{1}(V) \rtimes C_{2}, \zit/2) \to H^{1} (\pi_{1}(V), \zit/2) \to  H^{2} (C_{2}, \zit/2) \to H^{2} (\pi_{1}(V) \rtimes C_{2}, \zit/2)\]
Since the projection $\pi_{1}(V) \rtimes C_{2} \to C_{2}$ is tautologically split-injective, it follows that the rightmost arrow is injective, from which we deduce that the leftmost arrow is surjective, as desired.
\end{proof}
\end{Lem}
\begin{Rem}\label{ins} The assumption that (\ref{want-split}) be split on the right is necessary in the above lemma; indeed, if one takes $V$ to be an elliptic curve and $\iota: V \to V$ to be translation by a $2$-torsion element, then we do not have the surjectivity of the pull-back $H^{1} (U, \zit/2) \to H^{1} (V, \zit/2)$ in this case (since $U$ is also an elliptic curve and the kernel is non-trivial).
\end{Rem}
\noindent Remark \ref{ins} suggests that (\ref{want-split}) will be split on the right if $V$ admits a compactification to which the action of $\iota$ extends (and is no longer fixed-point-free). This is essentially true; one instance in which Assumption \ref{ass} holds is the following. Suppose that $V$ admits a smooth compactification $Y$ for which $\iota$ extends to an action $Y$ and for which the complement $S:= Y \setminus V$ satifies:
\begin{enumerate}[label=(\alph*)]
\item\label{smooth} $S$ is a smooth (non-empty) closed subvariety of codimension $n>1$
\item\label{fix} every $p \in S$ is fixed by $\iota$ and the codifferential map acting on the cotangent spaces $\delta\iota_{p}: T_{p}^{*}Y \to T_{p}^{*}Y$ acts by $-1$ on $\ker \{T_{p}^{*}Y \xrightarrow{res}  T_{p}^{*}D\}$
\end{enumerate}
We forget the variety structures and view $S$, $Y$, and $V$ as real manifolds. Then, by standard results in topology, there exists a tubular neighborhood $B$ of $S$ in $Y$ (the usual definition; c.f., \cite{BT} p. 65-66); by shrinking $B$ if necessary, we may assume that the fibers of the tautological map $B \to S$ are real Euclidean-open $2n$-balls (with spherical boundary), that $B$ is stable under the action of $\iota$ and that $\iota$ acts on the fibers of $B \to S$ by $-1$ (since we have \ref{fix} by assumption). By shrinking again if necessary, we can assume that $\partial B$ is stable under the action of $\iota$. The fibers of $\partial B \to S$ are $(2n-1)$-spheres and since $\iota$ acts on the fibers of $B \to S$ by $-1$, it follows that the quotient $P:= \partial B/\iota \to S$ is an $\mathbb{R}\mathbb{P}^{2n-1}$-fibration over $S$. Then, as in the first paragraph of the proof of Theorem 1 in \cite{Sp}, the sequence (\ref{want-split}) is split on the right in this case. Indeed, we have a commutative diagram with rows exact:
\begin{equation}\begin{tikzcd}\label{obv-cd} 
1 \arrow{r} & \pi_{1} (B) \arrow{r} \arrow{d} &  \pi_{1} (P) \arrow{r} \arrow{d} &  C_{2} \arrow{r} \arrow{d}{=} & 1\\
1 \arrow{r} & \pi_{1} (V) \arrow{r} &  \pi_{1} (U) \arrow{r}  &  C_{2} \arrow{r}  & 1 \end{tikzcd} \end{equation}
where the vertical arrows are push-forward maps. So, one is reduced to showing that the top sequence is split on the right, but this follows from the fact that $P \to S$ is an $\mathbb{R}\mathbb{P}^{2n-1}$-fibration (and the long exact sequence of homotopy groups of a fibration); note that $n>1$, so the fibers have fundamental groups $\cong C_{2}$.\\
\indent As the above discussion shows, one encounters the \'etale descent issue even in degree $1$. To avoid this issue in the application of Lemma \ref{litmus}, we restrict our attention to cases in which the answer is ``yes" trivially. In particular, we have
\begin{Cor}\label{triv} If $\iota$ acts trivially on $H^{1} (V, \zit/2)$, Assumption \ref{ass} holds and the cup product
\[ \wedge^{i} H^{1} (V, \zit/2) \to H^{i} (V, \zit/2) \]
is surjective, then the pull-back map $H^{i} (U, \zit/2) \xrightarrow{\pi^{*}} H^{i} (V, \zit/2)$ is surjective.
\begin{proof} This follows directly from Lemma \ref{lem-1}.
\end{proof}
\end{Cor}
\begin{Cor}\label{Span} Let $A$ be an Abelian variety over $\cit$ of dimension $g$ with a non-trivial involution $\iota: A\to A$ and let $\mathring{A} \subset A$ be the maximal Zariski open subset on which $\iota$ acts freely; denote the corresponding quotient $U := \mathring{A}/\iota$. Further, suppose that $\iota$ acts trivially on $H^{1} (A, \zit/2)$ and that $S= A \setminus \mathring{A}$ is a smooth (non-empty) closed subvariety of $A$ of codimension $\geq j >1$ for which the codifferential $\delta\iota$ satisfies condition \ref{fix} on the previous page, then the pull-back map:
\[ H^{i} (U, \zit/2) \xrightarrow{\pi^{*}} H^{i} (\mathring{A}, \zit/2) \cong H^{i} (A, \zit/2) \]
is surjective for $i <2j-1$.
\begin{proof} Using the Gysin sequence and the fact that the codimension of $A \setminus \mathring{A}$ in $A$ is $\geq j >1$, we deduce that the restriction 
\[ \wedge^{i} H^{1} (A, \zit/2) \xrightarrow[\cong]{\cup} H^{i} (A, \zit/2) \to H^{i} (\mathring{A}, \zit/2) \] 
is an isomorphism for $i < 2j-1$. In particular, the cup product $\wedge^{i} H^{1} (\mathring{A}, \zit/2) \to H^{i} (\mathring{A}, \zit/2)$ is surjective in these degrees. Moreover, by the paragraph following Lemma \ref{lem-1}, Assumption \ref{ass} applies in this case. The result now follows by Corollary \ref{triv}.
\end{proof}
\end{Cor}
\begin{Rem}\label{3-part} 
The main challenge in applying the above descent strategy more generally is in ensuring that Problem \ref{hard} admits a positive solution. The obstruction is given by the non-vanishing of the push-forward. Indeed, there is a long exact sequence:
\[ \ldots \to H^{i} (U, \zit/2) \xrightarrow{\pi^{*}} H^{i} (V, \zit/2) \xrightarrow{\pi_{*}} H^{i} (U, \zit/2) \to \ldots \]
arising from the short exact sequence of $\pi_{1}(U)$-modules:
\[ 0 \to \zit/2 \to \zit/2[C_{2}] \to \zit/2 \to 0\]
where the rightmost non-zero arrow is the trace map. The above long exact sequence has been used by many authors (see, for instance, \cite{KPS}, \cite{Sk} \S 2). Then, Problem \ref{hard} is equivalent to asking if $\pi_{*}$ vanishes on $H^{i} (V, \zit/2)$ whenever $\iota$ acts trivially on $H^{i} (V, \zit/2)$.
\end{Rem}
\noindent We will also need the extension lemma below that shows that one can promote unramified cohomology cycles from a smooth variety to its compactification under favorable circumstances.
\begin{Lem}\label{ext} Suppose that $U$ is a smooth irreducible variety over $\cit$ and for some integer $i$ for which the natural map 
\begin{equation} H^{i} (U, \zit/2(m)) \to H^{i}_{nr} (U, \zit/2(m))\label{image} \end{equation}
is non-zero. Suppose there exists some smooth projective compactification $X$ for which $D= X \setminus U$ is a divisor on $X$ with smooth irreducible components $D_{k}$ such that $H^{i-1}_{nr} (D_{k},\zit/2(m-1))= 0$ for all $k$. Then, any $\alpha$ in the image of (\ref{image}) lies in $H^{i}_{nr} (X, \zit/2(m)) \subset H^{i}_{nr} (U, \zit/2(m)) $.
\begin{proof}
It will suffice to show that $\alpha$ extends to a class in $H^{i}_{nr} (X, \zit/2(m))$, to which end we have the following diagram:
\[\begin{tikzcd} H^{i} (U, \zit/2(m)) \arrow{r} \arrow{d} & H^{i}_{nr} (U, \zit/2(m)) \arrow{d}\\
\bigoplus_{k} H^{i-1} (D_{k}, \zit/2(m-1)) \arrow{r} & \bigoplus_{k} H^{i-1} (\cit(D_{k}), \zit/2(m-1))
\end{tikzcd}
\]
where the horizontal arrows are restriction maps and the vertical arrows are the sums of residue maps. Now, we observe that the image of the lower horizontal arrow is precisely $H^{i-1}_{nr} (D_{i}, \zit/2(m-1))$, which vanishes by assumption. Thus, it follows that the lower horizontal arrow is $0$. Hence, $\alpha$ lies in the kernel of the right vertical arrow, which is precisely $H^{i}_{nr} (X, \zit/2(m))$ (using (\ref{kernel})). This gives the desired result.
\end{proof}
\end{Lem}
\begin{Rem} Note that it is entirely possible for $H^{i} (X, \zit/2(m)) = 0$ even though $H^{i}_{nr} (X, \zit/2(m)) \neq 0$ (c.f., Remark \ref{Span-rem}).
\end{Rem}
\subsection{Kummer varieties}

\noindent \noindent We now apply the above results to the case of a Kummer variety. Suppose that $A$ is a complex Abelian variety of dimension $g$ that is defined over a number field with good ordinary reduction at a prime dividing $2$. Let $\iota: A \to A$ be the inversion map on $A$ and $\mathring{A} = A \setminus A[2]$. As above, set $U= \mathring{A}/\iota$. Then, we consider the corresponding desingularized Kummer variety $X$ of $A$; this is a minimal compactification of $U$ obtained in the usual way: first blow up $A[2]$ on $A$ and then take the quotient of the corresponding blow-up by the induced action of $\iota$.  
\begin{Cor}\label{Cor-to-Sp} Suppose that $A$ is a complex Abelian variety of dimension $g$ that is defined over a number field with good ordinary reduction at a prime dividing $2$. Then, in the notation of the previous paragraph, the image of the natural map
\begin{equation} H^{i} (U, \zit/2(m)) \to H^{i}_{nr} (U, \zit/2(m))\label{target}\end{equation}
is non-zero for $1 \leq i \leq g$.
\begin{proof} It suffices to do this for $m=0$. Since $\iota$ acts by $-1$ on $H^{1} (A, \zit)$, it acts trivially on $H^{1} (A, \zit/2)$ and the co-differential of $\iota$ satisfies condition \ref{fix} above. Corollary \ref{Span} then implies that the pull-back $H^{i} (U, \zit/2) \to H^{i} (\mathring{A}, \zit/2)$ is surjective for $1 \leq i < 2g-1$. (We note that this is also obtained on p. 158 of \cite{Sp}.) Moreover, by Corollary \ref{Cor-to-BE} 
\[  H^{i} (\mathring{A}, \zit/2(m)) \to H^{i}_{nr} (\mathring{A}, \zit/2(m))\]
is non-zero for $1 \leq i \leq g$. By Lemma \ref{litmus}, it follows that (\ref{target}) is also non-zero.
\end{proof}
\end{Cor}

\begin{Cor}\label{kum-cor} With the assumptions and notations above, suppose now that $g >1$. Then, for $2 \leq i \leq g$
\[ H^{i}_{nr} (X, \zit/2(m)) \neq 0 \]  
\begin{proof} By Corollary \ref{Cor-to-Sp}, the natural map
\[ H^{i} (U, \zit/2(m)) \to H^{i}_{nr} (U, \zit/2(m))\]
is non-zero for $1 \leq i \leq g$. Hence by Lemma \ref{ext}, it suffices to check that $H^{i-1}_{nr} (D_{k}, \zit/2(m-1)) = 0$ for $2 \leq i \leq g$ for all irreducible components $D_{k}$ of $X \setminus U$. However, this is trivially true; indeed, $D_{k} \cong \pit^{g-1}$ for all $k$, and we certainly have $H^{i-1}_{nr} (\pit^{g-1}, \zit/2(m-1)) = 0$ for $2 \leq i \leq g$.
\end{proof}
\end{Cor}
\begin{Rem}\label{Span-rem}
It ought to be noted that the main result of \cite{Sp} is that $H^{odd} (X, \zit)$ vanishes and that $H^{*} (X, \zit)$ is torsion-free. In particular, $H^{i} (X, \zit/2)$ vanishes when $i$ is odd even though $H^{i}_{nr} (X, \zit/2)$ does not.
\end{Rem}
\section{Application to the integral Hodge conjecture}
Suppose now that $X$ is a smooth and projective variety over $\cit$. Moreover, let $H^{2m} (X, \zit(m))_{alg}$ denote the image of the degree $m$ cycle class map in $H^{2m} (X, \zit(m))$. Then, consider the group
\[ Z^{2m} (X) = H^{2m} (X, \zit(m))/H^{2m} (X, \zit(m))_{alg}\]
Observe that the $n$-torsion subgroup $Z^{2m} (X)[n]$ is generated by the degree $2m$ cohomology classes $\alpha$ for which some multiple is algebraic. In particular, any such $\alpha$ is of type $(m,m)$, and so whenever $Z^{2m} (X)[n] \neq 0$, the integral Hodge conjecture fails. Then, we have the following result relating degree $3$ unramified cohomology to the failure of the integral Hodge conjecture in degree $4$:
\begin{Thm}[Colliot-Th\'el\`ene and Voisin, \cite{CTV} Th\'eor\`eme 3.7]\label{CTV} Let $X$ be a smooth projective variety over $\cit$. Then, for every integer $n$, there is a short exact sequence:
\[ 0 \to  H^{3}_{nr} (X, \zit(2)) \otimes \zit/n \to H^{3}_{nr} (X, \zit/n(2)) \to Z^{4} (X)[n] \to 0\]
\end{Thm}
\begin{Rem}\label{BS} In light of this result, a strategy for finding a counterexample would be to find a smooth projective variety $X$ for which $H^{3}_{nr} (X, \zit(2)) = 0$ but for which $H^{3}_{nr} (X, \zit/n(2)) \neq 0$ for some $n$. The condition that $H^{3}_{nr} (X, \zit(2)) = 0$ is satisfied if, for instance, $CH_{0} (X \setminus S)_{\qit}=0$ for some closed surface $S \subset X$ (as noted in \cite{CTV}). This follows from the proof of Theorem 1 in \cite{BS}.
\end{Rem}
\noindent We now arrive at the following consequence of the above, which proves Theorem \ref{main}.
\begin{Cor}\label{main-cor} There exist complex Kummer varieties of dimension $3$ defined over number fields possessing non-algebraic non-torsion $(2,2)$ cohomology classes. 
\begin{proof} Let $A$ be an Abelian variety of dimension $3$ that is defined over a number field and that has good reduction at some prime dividing $2$. Let $X$ be the corresponding Kummer variety. Now, we have that $CH_{0} (X \setminus S)_{\qit}=0$ for some closed surface $S \subset X$ by \cite{BS} \S 4 Example (1) (see also \cite{B} Prop. 7). Thus, by Remark \ref{BS}, we deduce that $H^{3}_{nr} (X, \zit(2)) = 0$. So, by Corollary \ref{kum-cor} and Theorem \ref{CTV}, it then follows that there exists some nontrivial $\gamma \in Z^{4} (X)[2]$. That $\gamma$ is non-torsion follows from Remark \ref{Span-rem}.
\end{proof} 
\end{Cor}
\begin{Rem} As noted in Remark \ref{many}, there is a Zariski dense subset of the moduli space of ppav's of dimension $3$ consisting of ppav's that are defined over $\qit$ with good ordinary reduction. In particular, the corresponding Kummer varieties fail the integral Hodge conjecture.
\end{Rem}
\noindent Retaining the notation of the previous result, it does not seem clear which algebraic cycle on $X$ is responsible for the failure of the integral Hodge conjecture; i.e., which $\alpha \in H^{2,2} (X, \zit)$ is non-algebraic. We can, however, give the following heuristic that suggests that the counterexample arises from the exceptional locus of $X$. Indeed, let $D := X \setminus U \xhookrightarrow{i} X$ and consider the Gysin sequence:
\[ 0 = H^{3} (X, \zit(2)) \to H^{3} (U, \zit(2)) \to H^{2} (D, \zit(1)) \xrightarrow{i_{*}} H^{4} (X, \zit(2)) \to H^{4} (U, \zit(2)) \]
Since $H^{2} (D, \zit(1))$ is torsion-free, it is easy to see that $i_{*}$ in the above sequence is injective (since the analogous map $\otimes \qit$ is injective), from which we deduce that $H^{3} (U, \zit(2))$. Thus, the Bockstein sequence gives an isomorphism
\[H^{3} (U, \zit/2(2)) \cong  H^{4} (U, \zit(2))[2] \]
In particular, there is some $\gamma \neq 0 \in H^{4} (U, \zit(2))[2]$ corresponding to any $\tilde{\gamma} \in H^{3} (U, \zit/2(2))$ whose image under
\begin{equation} H^{3} (U, \zit/2(2)) \to H^{3} (\mathring{A}, \zit/2(2))  \end{equation}
is non-zero (see Corollary \ref{Span}), and by the proof of Corollary \ref{Cor-to-Sp}, this gives an unramified cycle. Since $H^{3} (D, \zit) = 0$, the Gysin sequence above implies that the restriction map $H^{4} (X, \zit(2)) \to H^{4} (U, \zit(2))$ is surjective. So, we let $\alpha$ be a lift of $\gamma$ to $H^{4} (X, \zit(2))$. Since $\gamma$ is $2$-torsion, there is some $\beta \in H^{2} (D, \zit(1))$ such that
\[ i_{*}\beta = 2\cdot \alpha \in H^{4} (X, \zit(2)) \]
The cycle $\beta$ is certainly algebraic since all the components of $D$ are algebraic since all the components of $D$ are isomorphic to $\pit^{2}$. The speculation then seems to be $\alpha$ is not algebraic, but it is not clear how to prove this directly.\\
\indent As a matter of completeness, we mention the following straightforward result that shows that the failure of the integral Hodge conjecture in the Kummer case can only arise from $2$-torsion in $Z^{4} (X)$:
\begin{Lem} For any complex Kummer variety $X$ of dimension $3$, $Z^{4} (X)[n] = 0= H^{3}_{nr} (X, \zit/n(2))$ for all $n$ odd. 
\begin{proof} Suppose $\alpha \in H^{2,2} (X, \zit)$ is such that $n\cdot \alpha \in H^{4} (X, \zit(2))$ is algebraic. Let $A$ be the corresponding Abelian variety, $\tilde{A}$ its blow-up along $A[2]$ and $\pi: \tilde{A} \to A$ the corresponding quotient map. By \cite{G} Chapter 2 the integral Hodge conjecture holds for $A$ (and, hence, also for $\tilde{A}$), so we deduce that $\pi^{*}\alpha \in H^{2,2} (\tilde{A}, \zit)$ is algebraic. It follows that
\[ 2\cdot \alpha = \pi_{*}\pi^{*}\alpha \in H^{2,2} (X, \zit) \]
is also algebraic. Since $n\cdot \alpha$ is also algebraic and $n$ is odd, this implies that $\alpha \in H^{2,2} (X)$ is algebraic. Hence, $Z^{4} (X)[n] = 0$. Using the fact that $H^{3}_{nr} (X, \zit(2))$ vanishes by the proof of Corollary \ref{main-cor}, it follows from Theorem \ref{CTV} that the group $H^{3}_{nr} (X, \zit/n(2))$ vanishes.
\end{proof}
\end{Lem}
\section{A higher-dimensional generalization}
\noindent It is tempting to think that Corollary \ref{main-cor} might also hold if $g>3$. This should be true, at least conjecturally. Indeed, one needs to prove that $H^{3}_{nr} (X, \zit(2))=0$; however, as in the proof of \cite{BS} Theorem 1, it suffices to show that $H^{3}_{nr} (X, \qit(2))=0$. This is true, provided that one believes outstanding conjectures about cycles on Abelian varieties, which we describe below. For an Abelian variety $A$ (over $\cit$) of dimension $g$, there is the Beauville decomposition \cite{B} of the Chow group:
\[ CH^{*} (A)_{\qit} = \bigoplus_{s, i} CH^{i}_{(s)} (A)_{\qit}\]
where $CH^{i}_{(s)} (A)_{\qit} := \{ \alpha \in CH^{i} (A)_{\qit} \ | \ \ n^{*}\alpha = n^{2i-s}\alpha \ \forall n \in \zit \}$.
Here $n^{*}$ denotes the induced action of multiplication by $n$ on $A$. It is known that $CH^{i}_{s} (A)_{\qit}$ vanishes unless $i-g< s <i$ (see first result of op. cit.). There is the following expectation:
\begin{Conj}[Beauville] $CH^{i}_{(s)} (A)_{\qit}$ vanishes for $s <0$ and the cycle class map
\[ CH^{i}_{(0)} (A)_{\qit} \to H^{2i} (A, \qit(i)) \] 
is injective for all $i$.
\end{Conj}
\noindent Since $n^{*}$ acts by $n^{2i}$ on $H^{2i} (A, \qit(i))$, the cycle class map vanishes on the summands $CH^{i}_{(s)} (A)_{\qit}$ for $s \neq 0$. So, it is necessary to restrict to the summand $CH^{i}_{(0)} (A)_{\qit}$ to have injectivity. The first statement of Beauville's conjecture is known for $i=0,1, g-2,g-1, g$ (by Prop. $3$ of \cite{B}), while the second statement is known for $i=0,1,g-1, g$. Apart from this, the conjecture remains open in general.\\
\indent As to the other extreme, there is the following basic lemma that is surely well-known to the experts but for which a reference was not found:
\begin{Lem}\label{alg} For all $1 \leq i \leq g$, $CH^{i}_{(i)} (A)_{\qit} \subset CH^{i}_{alg} (A)_{\qit}$, where $CH^{i}_{alg} (A)_{\qit}$ is the subspace of algebraically-equivalent-to-zero cycles in $CH^{i} (A)_{\qit}$.
\begin{proof} This is certainly true for $i=1$, so we assume that $i \geq 2$, and this case then becomes an exercise in applying properties of the Fourier transform. To this end, let $\mathcal{P}$ be the Poincar\'e invertible sheaf on $A \times \hat{A}$, where $\hat{A} = Pic^{0} (A)$ is the dual Abelian variety. Also, let 
\[ \mathcal{F}: CH^{*} (A)_{\qit} \to CH^{*} (\hat{A})_{\qit}, \ \alpha \mapsto \pi_{\hat{A}*}(\pi_{A}^{*}\alpha\cdot ch(\mathcal{P})) \]
be the corresponding Fourier transform, where $ch(-)$ denotes the Chern character. By Prop. 2 of \cite{B}, we have that
\begin{equation} \mathcal{F}_{A} : CH^{i}_{(s)} (A)_{\qit} \xrightarrow{\cong} CH^{g-i+s}_{(s)} (A)_{\qit}\label{fourier} \end{equation}
In particular, we have that $\mathcal{F}_{A}(CH^{i}_{(i)} (A)_{\qit}) = CH^{g}_{(i)} (A)_{\qit}$. Using (\ref{fourier}), it also follows that $CH^{g}_{(0)} (A)_{\qit} = \qit\cdot[e]$, where $e \in A$ is the identity and that
\[CH^{1}_{(1)} (A)_{\qit} \xrightarrow{\cong} CH^{g}_{(1)} (A)_{\qit}  \]
which implies that $CH^{g}_{(1)} (A)_{\qit} \cong A \otimes \qit$ via the Albanese map. Hence, the Albanese kernel $I$ is precisely $\bigoplus_{i\geq 2} CH^{g}_{(i)}(A)_{\qit}$, and since $I \subset CH^{g}_{alg} (A)_{\qit}$, we deduce that $CH^{g}_{(i)}(A)_{\qit} \subset CH^{g}_{alg} (A)_{\qit}$ for all $i\geq 2$. Since $\mathcal{F}_{\hat{A}}\circ\mathcal{F}_{A} = (-1)^{g}$ and since 
\[ \mathcal{F}_{\hat{A}}(CH^{*}_{alg} (\hat{A})_{\qit}) \subset CH^{*}_{alg} (A)_{\qit} \]
it follows that $CH^{i}_{(i)} (A)_{\qit} \subset CH^{*}_{alg} (A)_{\qit}$.
\end{proof}
\end{Lem}
\noindent As usual, we let $CH^{i}_{hom} (A)_{\qit}$ denote the subspace of $CH^{i} (A)_{\qit}$ of nul-homologous cycles and the Griffiths group by $Griff^{i} (A)_{\qit} := CH^{i}_{hom} (A)_{\qit}/CH^{i}_{alg} (A)_{\qit}$.
\begin{Cor}\label{sec-last} Assuming Beauville's conjecture for $i=2$, $\iota$ acts as $-1$ on $Griff^{2} (A)_{\qit}$. Moreover, $Griff^{2} (X)_{\qit} = 0$, where $X$ the associated Kummer variety of $A$.
\begin{proof} Under the assumption of Beauville's conjecture for $i=2$, Beauville's decomposition in codimension $2$ becomes:
\[ CH^{2} (A)_{\qit} = CH^{2}_{(0)} (A)_{\qit} \oplus CH^{2}_{(1)} (A)_{\qit} \oplus CH^{2}_{(2)} (A)_{\qit}\]
By Lemma \ref{alg}, we have $CH^{2}_{(2)} (A)_{\qit} \subset CH^{2}_{alg} (A)_{\qit}$. There is then an iduced decomposition on the Griffiths group:
\[ Griff^{2} (A)_{\qit} = Griff^{2}_{(0)} (A)_{\qit} \oplus Griff^{2}_{(1)} (A)_{\qit}\]
Now, by Beauville's conjecture, we also have $CH^{2}_{hom} (A)_{\qit} \cap CH^{2}_{0} (A)_{\qit} = 0$, from which we deduce that $Griff^{2} (A)_{\qit} = Griff^{2}_{(1)} (A)_{\qit}$. In particular, $(-1)^{*}$ acts as $(-1)$ on $Griff^{2} (A)_{\qit}$, as was to be shown. The proof of the second statement is then the same as that of Prop. 7 in \cite{B}.
\end{proof}
\end{Cor}
\noindent This last result is the necessary step to generalize Corollary \ref{main-cor} to higher dimension.
\begin{Cor}\label{last} Suppose that Beauville's conjecture for $i=2$. Let $X$ be a very general complex Kummer variety of dimension $g>2$. Then, there exist non-algebraic non-torsion $(2,2)$ cohomology classes on $X$. Moreover, one can find such $X$ defined over any given number field.
\begin{proof} The proof of Corollary \ref{main-cor} works mutatis mutandis except that one has to verify that $H^{3}_{nr} (X, \qit(2))$, as noted at the beginning of this section. To this end, we note that from the local-to-global spectral sequence:
\[ H^{p}_{Zar} (X, \mathcal{H}^{q}_{X} (\qit(n)) \Rightarrow H^{p+q} (X, \qit(n)) \]
one obtains the following short exact sequence:
\[ H^{3} (X, \qit(2)) \to H^{3}_{nr} (X, \qit(2)) \to Griff^{2} (X) \otimes \qit \]
By Remark \ref{Span-rem}, the first group vanishes, so one is reduced to showing that $Griff^{2} (X) \otimes \qit$, which follows from Corollary \ref{sec-last}.
\end{proof}
\end{Cor}

\Addresses
\end{document}